\documentclass{article}

\usepackage{amssymb}

\newcommand{\real}{{\mathbb R}}
\newcommand{\iin}{\!\in\!}
\newcommand{\Lipplus}{\mbox{\rm BLip}^+}
\newcommand{\Ubru}[1]{{\bf U_b}(\ru #1)}
\newcommand{\Measru}[1]{{\bf M}(\ru #1)}
\newcommand{\UMeasru}[1]{{\bf M_u}(\ru #1)}
\newcommand{\sect}[1]{\setminus_{#1}}
\newcommand{\bigsep}{\mbox{\Large $|$}}
\newcommand{\conv}{\star}
\newcommand{\ru}{r\,}
\newcommand{\wstar}{weak\raisebox{1mm}{$\ast$}\ }
\newcommand{\rtrans}[2]{\rho^{#1}(#2)}
\newcommand{\compln}[1]{\widehat{#1}}
\newcommand{\compactnru}[1]{\overline{\ru #1}}
\newcommand{\qed}{\hfill $\Box$\vspace{3mm}}
\newcommand{\rp}{{\cal RP}}
\newcommand{\orbit}[1]{\mbox{\rm orb}(#1)}
\newcommand{\clorbit}[1]{\overline{\mbox{\rm orb}(#1)}}
\newcommand{\card}[1]{| #1 |}
\newcommand{\bcard}{\eta^{\textstyle \sharp}}
\newcommand{\bccard}{\eta}    
\newcommand{\scrO}{{\cal O}}
\newcommand{\psm}{\Delta}
\newcommand{\cf}{\mbox{\rm cf}}

\newtheorem{theorem}{Theorem}[section]
\newtheorem{lemma}[theorem]{Lemma}
\newtheorem{corollary}[theorem]{Corollary}
\newtheorem{question}{Question}

\title{Ambitable topological groups}

\author{Jan Pachl  \\[10pt]
\emph{Fields Institute} \\
\emph{Toronto, Ontario, Canada}}
\date{May 25, 2009~(version 4)}

\begin{document}
\maketitle

\begin{abstract}
A topological group is said to be {\em ambitable\/} if each
uniformly bounded uniformly equicontinuous set of functions on the
group with its right uniformity is contained in an ambit.
For $n=0,1,2, \ldots$,
every locally $\aleph_n$-bounded topological group is either
precompact or ambitable. In the familiar semigroups constructed
over ambitable groups, topological centres have an effective
characterization.
\end{abstract}


\section{Overview}
    \label{section:overview}

A topological group $G$ may be naturally embedded in larger spaces,
algebraically and topologically.
Two such spaces of particular interest in abstract harmonic analysis are
\begin{itemize}
\item
the norm dual of the space of bounded right uniformly continuous functions on $G$,
denoted here $\Measru{G}$ (also known as LUC$(G)\raisebox{1mm}{$\ast$}$); and
\item
the uniform compactification of $G$ with its right uniformity,
denoted here $ \compactnru{G} $ (also known as $ G^{\mbox{\rm\scriptsize LUC}} $
or $ G^{\cal LC} $).
\end{itemize}
It is customary to study these ``right'' versions of the two spaces;
the properties of the corresponding ``left'' versions are obtained by symmetry.

Both $\Measru{G}$ and $ \compactnru{G} $ are right topological semigroups.
In investigating their structure it is very helpful to have
a tractable characterization of their topological centres.
Feasible candidates for such characterizations are the space of uniform measures
$ \UMeasru{G} $ and the completion $ \compln{\ru G} $ of the right uniformity on $G$.

When $G$ is locally compact, $ \UMeasru{G} $ is the space of finite Radon measures on $G$,
and $ \compln{\ru G} $ is $G$ itself.
In this case,
Lau~\cite{Lau1986} and Lau and Pym~\cite{Lau-Pym1995} proved that
$ \UMeasru{G} $ and $ \compln{\ru G} = G $ are the topological centres of $\Measru{G}$
and $ \compactnru{G} $.
These characterizations generalized a number of previous results for special
classes of locally compact groups.

More recently, Neufang~\cite{Neufang2004} applied his factorization method to simplify the
proof of  Lau's result.
Then Ferri and Neufang~\cite{Ferri-Neufang2007} used a variant of the factorization method
to prove that $ \UMeasru{G} $ and $ \compln{\ru G} $ are the topological centres
of $\Measru{G}$ and $ \compactnru{G} $ for
$\aleph_0$-bounded (not necessarily locally compact) topological groups.

This paper deals with another variant of the factorization method, similar to that used by
Ferri and Neufang.
By definition, ambitable topological groups are those in which a suitable factorization
theorem holds;
equivalently, in the language of topological dynamics,
those in which every uniformly bounded uniformly equicontinuous
set of functions is contained in an ambit.
In such groups $ \UMeasru{G} $ and $ \compln{\ru G} $ are the topological centres
of $\Measru{G}$ and $ \compactnru{G} $.
Several classes of topological groups are shown to be ambitable.
In particular, if $n$ is a positive integer
then every locally $\aleph_n$-bounded group is either precompact or ambitable,
which yields a common generalization of the aforementioned results by Lau, Lau and Pym,
and Ferri and Neufang.

This is an updated preliminary version of the paper.
The final version will be published in Topology and its Applications in 2009.


\section{Basic definitions}
    \label{section:definitions}

All topological groups considered in this paper are assumed to be Hausdorff,
and all linear spaces to be over the field $\real$ of reals.

Most of the notation used here is defined in~\cite{Pachl2006}.
In particular,
if $p$ is a mapping from $ X \times Y $ to $ W $ then
$ \sect{x} p(x,y) $ is the mapping $ x \mapsto p(x,y) $ from $X$ to $W$
and $ \sect{y} p(x,y) $ is the mapping $ y \mapsto p(x,y) $ from $Y$ to $W$.

Let $G$ be a group, $f$ a real-valued function on $G$ and $ x \in G $.
Define
$\rtrans{x}{f}$ (the {\em right translation\/} of $f$ by $x$) to be the function
$ \sect{z} f(zx) $.
The set
$ \orbit{f} = \{ \rtrans{x}{f} \; | \; x \in G \} $
is the {\em (right) orbit\/} of $f$.
Denote by $ \clorbit{f} $
the closure of $ \orbit{f} $ in the product space $ \real^G $
(the set of real-valued functions on $G$ with the topology of pointwise convergence).

When $\psm$ is a pseudometric on $G$, define
\[
\Lipplus ( \psm )  =  \{ f: G \rightarrow \real \; \bigsep \; 0 \leq f(x) \leq 1 \;
\mbox{\rm and} \;
| f(x) - f(y) | \leq \psm(x,y) \; \mbox{\rm for all} \; x,y \in G  \} \;\; .
\]
Then $\Lipplus(\psm)$ is a compact subset of the product space $ \real^G $;
we always consider $\Lipplus(\psm)$
with this compact topology.

When $G$ is a topological group, denote by $ \rp (G) $ the set of all continuous
right-invariant pseudometrics on $G$.
The {\em right uniformity\/} on $G$ is the uniformity generated by $ \rp (G) $.
The set $G$ with the right uniformity is denoted $ \ru G $,
and the space of all bounded uniformly continuous real-valued functions on $ \ru G $
is denoted $ \Ubru{G} $.
The group $G$ is said to be {\em precompact\/} if the uniform space $ \ru G $ is precompact.

The following lemma summarizes several properties of $\rtrans{x}{f}$ needed
in this paper.
De Vries (\cite{deVries1993}, sec.~IV.5) provides a comprehensive treatment
of the role of $\rtrans{x}{f}$ in topological dynamics.

\begin{lemma}
Let $G$ be any topological group and $\psm \in \rp (G)$.
\begin{enumerate}
\item
If $f$ is a real-valued function on $G$ and $x,y \in $G then
$ \rtrans{xy}{f} = \rtrans{x}{\rtrans{y}{f}} $.
\item
If $f \in \Lipplus(\psm)$ and $x \in G$ then $\rtrans{x}{f} \in \Lipplus(\psm) $.
\item
The mapping $ (x,f) \mapsto \rtrans{x}{f} $ is continuous
from $ G \times \Lipplus(\psm) $ to $ \Lipplus(\psm) $.
\end{enumerate}
\end{lemma}

\noindent
{\bf Proof.}
1. $ \rtrans{xy}{f} (z) = f(zxy) = \rtrans{y}{f} (zx) = \rtrans{x}{\rtrans{y}{f}} (z) $.

2. If $ f \in \Lipplus(\psm) $ then
$ | \rtrans{x}{f} (z) - \rtrans{x}{f} (z') | = | f(zx) - f(z' x) |
\leq \psm(zx, z'x) = \psm(z,z' ) $.

3. To prove that the mapping $ (x,f) \mapsto \rtrans{x}{f} $ to $ \Lipplus(\psm) $
is continuous, it is sufficient to prove that the mapping
$ (x,f) \mapsto \rtrans{x}{f}(z) $ to $\real$ is continuous for each $ z \in G $.

Take any $ z \in G $, $ ( x_0 , f_0 ) \in G \times \Lipplus(\psm) $ and $ \varepsilon > 0 $.
The set
\[
U = \{ \; (x,f) \in G \times \Lipplus(\psm) \; \bigsep \; \psm (zx, z x_0 ) < \varepsilon
\;\; \mbox{\rm and} \;\;
| f(z x_0 ) - f_0 (z x_0 ) | < \varepsilon \; \}
\]
is a neighbourhood of $ ( x_0 , f_0 ) $ in $ G \times \Lipplus(\psm) $.
For $ (x,f) \in U $ we have
\[
| \rtrans{x}{f} (z) - \rtrans{x_0}{f_0} (z) |
= | f(zx) - f_0 (z x_0) |
\leq | f(zx) - f(z x_0 ) | + | f (z x_0 ) - f_0 ( z x_0 ) | < 2 \varepsilon .
\]
Thus the mapping $ (x,f) \mapsto \rtrans{x}{f}(z) $ is continuous at $ ( x_0 , f_0 ) $.
\qed

When $G$ and $\psm$ are as in the lemma,
$\Lipplus(\psm)$ with the action $ (x,f) \mapsto \rtrans{x}{f} $ is a compact $G$-flow,
in the terminology of topological dynamics~\cite{deVries1993}.
If $ f \in \Ubru{G} $ then there exists $ \psm \in \rp(G) $ such that
$ sf + t \in \Lipplus(\psm) $ for some $ s,t \in \real$.
Thus $ s \cdot \, \clorbit{f} + t = \clorbit{sf+t} \subseteq \Lipplus(\psm) $ and therefore
the set $ \clorbit{f} $ is compact in the topology of pointwise convergence,
and $\clorbit{f}$ with the action $ (x,f) \mapsto \rtrans{x}{f} $ is also a compact $G$-flow.

Recall that a compact $G$-flow is an {\em ambit\/}
if it contains an element with dense orbit (\cite{deVries1993}, IV.4.1).
For a fixed $G$,
all ambits can be constructed from those of the form $\clorbit{f}$,
where $f \in \Ubru{G}$ (\cite{deVries1993}, IV.5.8).

Say that a topological group $G$ is {\em ambitable\/} if every $\Lipplus(\psm)$,
where $ \psm \in \rp(G) $, is contained in an ambit within $\Ubru{G}$.
In other words, $G$ is ambitable iff for each $ \psm \in \rp(G) $
there exists $ f \in \Ubru{G} $ such that
$ \Lipplus(\psm) \subseteq \clorbit{f} $.

\begin{theorem}
    \label{theorem:precompact}
No precompact topological group is ambitable.
\end{theorem}

\noindent
{\bf Proof.}
Let $G$ be precompact.
Fix any $ f \iin \Ubru{G} $,
and define $ \psm ( x , x' ) = \displaystyle\sup_{y \in G } | f(xy) - f(x'y) | $
for $ x,x' \iin G$.
Then $ \psm \iin \rp(G) $, and
there is a finite set $ F \subseteq G $ such that $ \psm ( x , F ) < \frac{1}{3} $
for every $ x \iin G $.
Consider the constant functions $0$ and $1$.
If $ 0 \iin \clorbit{f} $ then
there is $ y \iin G $ such that
$ \rtrans{y}{f} (x) < \frac{1}{3} $ for every $ x \iin F $,
hence $ \rtrans{y}{f} (x) < \frac{2}{3} $ for every $ x \iin G $.
Thus $ f(x) < \frac{2}{3} $ for every $ x \iin G $, and $ 1 \not\in \clorbit{f} $.
This proves that there is no $ f \iin \Ubru{G} $ for which $ 0,1 \iin \clorbit{f} $.
\qed

\begin{question}
\label{question:ambitable}
Is every topological group either precompact or ambitable?
\end{question}

This question is motivated by investigations of topological centres in certain semigroups
that arise in functional analysis.
The connection is explained in section~\ref{section:centres} below.


\section{Cardinal functions}
    \label{section:cardfunctions}

The reader is referred to Jech~\cite{Jech2002} for definitions regarding cardinals.
The cardinality of a set $X$ is $\card{X}$.
The cardinal successor of a cardinal $\kappa$ is $\kappa^+$.
The smallest infinite cardinal is~$\aleph_0$, and
$ \aleph_{n+1} = \aleph_n^+ $.
The smallest cardinal larger than $\aleph_n$ for $n=0,1,2, \ldots $ is $\aleph_\omega$.

Let $G$ be a group and $\psm$ a pseudometric on $G$.
Sufficient conditions in the next section are expressed
in terms of three cardinal functions:
\newpage{
\begin{itemize}
\item
$d( \psm )$, the {\em $\psm$-density\/} of $G$
(the smallest cardinality of a $\psm$-dense subset of $G$);
\item
$ \bcard( \psm ) $,
the smallest cardinality of a set $ P \subseteq G $ such that
\[
G = \bigcup_{p \in P} \;
\{ \; x \in G \; | \; \psm ( p, x ) \leq 1 \; \} \; ;
\]
\item
$ \bccard( \psm ) $,
the smallest cardinality of a set $ P \subseteq G $ for which there exists
a finite set $ Q \subseteq G $ such that
\[
G = \bigcup_{q \in Q} \; \bigcup_{p \in P} \;
\{ \; x \in G \; | \; \psm ( p, qx ) \leq 1 \; \} \; .
\]
\end{itemize}
}

The following lemma collects basic facts about these three functions.
Proofs follow directly from the definition.

\begin{lemma}
    \label{lemma:cardftions}
Let $\psm$ be a pseudometric on a group $G$.
Let $ B = \{ x \iin G \; | \; \psm(e,x) \leq 1 \} $,
where $e$ is the identity element of $G$.
\begin{enumerate}
\item
    \label{lemma:cardftions:part1}
$ \bccard( \psm ) \leq \bcard( \psm ) \leq d( \psm ) $.
\item
    \label{lemma:cardftions:part2}
$ \displaystyle d( \psm ) = \lim_{k \rightarrow \infty} \bcard( k \psm ) $.
\item
If $\psm'$ is another pseudometric on $G$ such that $ \psm \leq \psm' $ then
$ \bccard(\psm) \!\leq \bccard(\psm') $,
$ \bcard(\psm) \!\leq \bcard(\psm') $
and $ d(\psm) \!\leq d(\psm') $.
\item
If $ \psm $ is left-invariant and $ \bccard ( \psm ) \geq \aleph_0 $ then
$ \bccard( \psm ) = \bcard( \psm ) $.
\item
    \label{lemma:cardftions:part5}
If $ \psm $ is right-invariant then $ \bcard( \psm ) $ is
the smallest cardinality of a set $ P \subseteq G $ such that $ G = B P $
and $ \bccard( \psm ) $ is
the smallest cardinality of a set $ P \subseteq G $ for which there exists
a finite set $Q$ such that $ G = Q B P $.
\end{enumerate}
\end{lemma}

It is easy to see that a topological group $G$ is precompact if and only if $ \bcard(\psm) $
is finite for each \mbox{$ \psm \iin \rp(G) $}.
Part~\ref{theorem:bctobcc:part1} in
Theorem~\ref{theorem:bctobcc} below yields a stronger statement:
$G$ is precompact if and only if $ \bccard(\psm) $
is finite for each \mbox{$ \psm \iin \rp(G) $}.
This is equivalent to the theorem of Uspenskij (\cite{Uspenskij2001}, p.~338;
\cite{Uspenskij2008}, p.~1581),
for which a simple proof was given by Bouziad and Troallic~\cite{Bouziad-Troallic2007}.
Ferri and Neufang~\cite{Ferri-Neufang2007} gave another proof
using a result of Protasov~(\cite{Protasov1997}, Th.~11.5.1).

The case of finite $P$ in the next Lemma is due to Bouziad and Troallic
(\cite{Bouziad-Troallic2007}, Lemma~4.1).
The proof below is a straightforward generalization
of their approach, which in turn was adapted from Neumann~\cite{Neumann1954}.

\begin{lemma}
   \label{lemma:partition}
Let $G$ be a group,
$ P \subseteq G $,
and $ A_k \subseteq G $ for $ 1 \leq k \leq n $.
If $ G = \bigcup_{k=1}^n A_k P $
then there are a set $ P' \subseteq G $ and $j$, $ 1 \leq j \leq n $,
such that $ G = A_j^{-1} A_j P' $ and
\begin{itemize}
\item
if $P$ is finite then so is $P'$, and
\item
if $P$ is infinite then $\card{P'} \leq \card{P} $.
\end{itemize}
\end{lemma}

\noindent
{\bf Proof} proceeds by induction in $n$.
When $n=1$, the statement is true with $j=1$ and $P' = P$.

For the induction step, let $m \geq 1$ and assume that the statement in the lemma
is true for $n = m$.
Let $ P \subseteq G $ and $ A_1 , A_2, \ldots , A_{m+1} \subseteq G $ be such that
$ G = \bigcup_{k=1}^{m+1} A_k P $.

If $ G = A_{m+1}^{-1} A_{m+1} P $ then set $ j = m+1 $ and $ P' = P $.

On the other hand, if $ G \neq A_{m+1}^{-1} A_{m+1} P $ then take any
$ x \iin G \setminus A_{m+1}^{-1} A_{m+1} P $.
It follows that $ A_{m+1} x \cap A_{m+1} P = \emptyset $,
and $ A_{m+1} \subseteq \bigcup_{k=1}^m A_k P x^{-1} $.
Thus
\[
G = \bigcup_{k=1}^m A_k ( P \, \cup \, P x^{-1} P ) \; .
\]
By the induction hypothesis, there are $ P' \subseteq G $ and $j$, $ 1 \leq j \leq m $,
such that $ G = A_j^{-1} A_j P' $, $P'$ is finite if $P$ is, and
$\card{P'} \leq \card{P} $ if $P$ is infinite.

Thus in either case the statement holds for $n=m+1$.
\qed

\begin{theorem}
   \label{theorem:bctobcc}
Let $G$ be a topological group, and $\psm \iin \rp(G)$.
\begin{enumerate}
\item
   \label{theorem:bctobcc:part1}
If $\bccard(\psm) $ is finite then $\bcard(\frac{1}{2} \psm)$ is finite.
\item
    \label{theorem:bctobcc:part2}
If $\bccard(\psm) $ is infinite then $\bcard(\frac{1}{2} \psm) \leq \bccard(\psm)$.
\end{enumerate}
\end{theorem}

\noindent
{\bf Proof.}
Let $ B = \{ x \iin G \; | \; \psm(e,x) \leq 1 \} $, where $e$ is the identity element of $G$.
If $ y,z \iin B $ then $ y^{-1}z \in \{ \;x \iin G \; | \; \frac{1}{2} \psm (e,x) \leq 1 \; \} $,
because
\[
\psm( e, y^{-1}z ) = \psm( z^{-1}, y^{-1} ) \leq \psm( z^{-1}, e ) + \psm( y^{-1}, e )
= \psm(e,z) + \psm(e,y) \leq 2 .
\]

By part~\ref{lemma:cardftions:part5} of Lemma~\ref{lemma:cardftions},
there are sets $ P, Q \subseteq G $ such that
$Q$ is finite, $ \card{P} = \bccard(\psm) $
and $ G = QBP $.
By Lemma~\ref{lemma:partition}, there are $ q \iin Q $ and $ P' \subseteq G $ such that
$P'$ is finite if $P$ is,
$\card{P'} \leq \card{P} $ if $P$ is infinite, and
\[
G = ( q B )^{-1} q B P' = B^{-1} B P' \subseteq \{ \;x \iin G \; | \;
{\textstyle \frac{1}{2}} \psm (e,x) \leq 1 \; \} \; P'
\]
which shows that $\bcard(\frac{1}{2} \psm) \leq | P' |$.
If $\bccard(\psm) $ is finite then $|P|$ and $|P'|$ are finite
and therefore $\bcard(\frac{1}{2} \psm)$ is finite.
If $\bccard(\psm) $ is infinite then
$\bcard(\frac{1}{2} \psm) \leq | P' | \leq |P| = \bccard(\psm)$.
\qed

\begin{corollary}
   \label{corollary:bctobcc}
Let $G$ be a topological group, and $\psm \iin \rp(G)$.
If $ d( \psm ) \! > \aleph_0 $ then
\[
\displaystyle d( \psm ) = \lim_{k \rightarrow \infty} \bccard( k \psm ) \; .
\]
\end{corollary}

\noindent
{\bf Proof.}
Combine parts~\ref{lemma:cardftions:part1} and~\ref{lemma:cardftions:part2} of
Lemma~\ref{lemma:cardftions} with
part~\ref{theorem:bctobcc:part2} of Theorem~\ref{theorem:bctobcc}.
\qed

Let $\kappa$ be an infinite cardinal.
Following the terminology of Guran~\cite{Guran1981}
(see also sec.~9 in~\cite{Arkhangelskii1981}),
say that a topological group $G$ is {\em $\kappa$-bounded\/}
if for every neighbourhood $U$ of the identity element in $G$ there
exists a set $H \subseteq G $ such that $ \card{H} \leq \kappa $ and $ U H = G $.

{\samepage
\begin{lemma}
    \label{lemma:bounded}
Let $\kappa$ be an infinite cardinal.
The following conditions for a topological group $G$ are equivalent:
\begin{description}
\item[{\it (i)}]
$G$ is $\kappa$-bounded;
\item[{\it (ii)}]
$ d(\psm) \leq \kappa $ for every $ \psm \iin \rp(G) $;
\item[{\it (iii)}]
$ \bcard(\psm) \leq \kappa $ for every $ \psm \iin \rp(G) $;
\item[{\it (iv)}]
$ \bccard(\psm) \leq \kappa $ for every $ \psm \iin \rp(G) $.
\end{description}
\end{lemma}
}  

\noindent
{\bf Proof.}
The family of all sets
$ \{ x \iin G \; | \; \psm(e,x) \leq 1 \} $,
where $ \psm \iin \rp(G)$,
is a basis of neighbourhoods of the identity element $e$ in $G$.
Therefore (i) $\Leftrightarrow$ (iii), by part~\ref{lemma:cardftions:part5}
in Lemma~\ref{lemma:cardftions}.

(ii) $\Rightarrow$ (iii) $\Rightarrow$ (iv) by part~\ref{lemma:cardftions:part1}
in Lemma~\ref{lemma:cardftions},
and (iv) $\Rightarrow$ (ii) by Corollary~\ref{corollary:bctobcc}.
\qed

Say that a topological group $G$ is {\em locally $\kappa$-bounded\/} if its identity
element has a neighbourhood
$U$ such that for each $ \psm \iin \rp(G) $ there is a
$\psm$-dense subset $H$ of $U$, $ \card{H} \leq \kappa$.
When $\kappa$ is an infinite cardinal,
every locally compact group is locally $\kappa$-bounded,
and so is every $\kappa$-bounded group.


\section{Sufficient conditions}
    \label{section:answers}

This section contains several sufficient conditions for a topological group
to be ambitable.
For each such condition we prove a slightly stronger property;
namely, that for every $ \psm \iin \rp (G) $ there exists
$ \psm' \iin \rp (G) $, $ \psm' \geq \psm $ such that
$ \Lipplus ( \psm' ) $ is an ambit.
The key result is Lemma~\ref{lemma:density},
which is another form
of the factorization theorems of Neufang~\cite{Neufang2004} and
Ferri and Neufang~\cite{Ferri-Neufang2007}.

\begin{lemma}
    \label{lemma:finsets}
Let $\psm$ be a pseudometric on a group $G$ such that $ \bccard(\psm) \geq \aleph_0 $.
Let $ A $ be a set of cardinality $ \bccard(\psm) $,
and for each $ \alpha \iin A $ let $ F_\alpha $ be a non-empty finite subset of $G$.
Then there exist elements $ x_\alpha \iin G $ for $ \alpha \iin A $ such that
$ \psm ( F_\alpha x_\alpha , F_\beta x_\beta ) > 1 $ whenever
$ \alpha , \beta \iin A $, $ \alpha \neq \beta $.
\end{lemma}

\noindent
{\bf Proof.}
Without loss of generality, assume that $A$ is the set of ordinals smaller than
the first ordinal of cardinality $ \bccard(\psm) $.
The construction of $ x_\alpha $ proceeds by transfinite induction.
For $ \gamma \in A $, let $ S( \gamma ) $ be the statement
``there exist elements $ x_\alpha \in G $ for all $ \alpha \leq \gamma $ such that
$ \psm ( F_\alpha x_\alpha , F_\beta x_\beta ) > 1 $ whenever
$ \alpha < \beta \leq \gamma \;$.''

Any choice of $ x_0 \in G $ makes $S(0)$ true.
Now assume that $\gamma \in A$, $\gamma > 0$,
and $S(\gamma')$ is true for all $ \gamma' < \gamma $.
We want to prove $S(\gamma)$.

Since $ \bccard(\psm) \geq \aleph_0 $
and the cardinality of $\gamma$ is less than $ \bccard (\psm) $,
from the definition of $ \bccard(\psm) $ we get
\[
G \neq  \bigcup_{q \in F_\gamma} \;
        \bigcup_{\alpha < \gamma} \;
        \bigcup_{p \in F_\alpha x_\alpha} \;
        \{ \; x \in G \; | \; \psm ( p, qx ) \leq 1 \; \} \; .
\]
Thus there exists $ x_\gamma \iin G $ such that $ \psm(p, q x_\gamma ) > 1 $ for
all $ q \in F_\gamma $ and all $ p \in F_\alpha x_\alpha $ where $ \alpha < \gamma $.
That means $ \psm ( F_\alpha x_\alpha , F_\gamma x_\gamma ) > 1 $
for all $ \alpha < \gamma $.
\qed

\begin{lemma}
    \label{lemma:largeorbit}
Let $G$ be a topological group, $ \psm \iin \rp(G)$ and
$ \bccard(\psm) \!\geq \aleph_0 $.
If $\scrO$ is a collection of non-empty open subsets of $ \Lipplus(\psm) $
and $ \card{\scrO} \leq \bccard (\psm) $,
then there exists $ f \iin \Lipplus(\psm) $
such that $ \orbit{f}$ intersects every set in $\scrO$.
\end{lemma}

\noindent
{\bf Proof.}
Without loss of generality,
assume that every set in $\scrO$ is a basic neighbourhood.
Thus each $U \in \scrO$ is of the form
\[
U = \{ \; f \in \Lipplus (\psm) \; \bigsep \;\; | f(x) - h_U (x) | < \varepsilon_U
 \; \mbox{\rm for} \; x \in F_U \; \}
\]
where $ F_U \subseteq G $ is a finite set, $h_U \iin \Lipplus(\psm) $, and
$ \varepsilon_U > 0 $.

By Lemma~\ref{lemma:finsets} with $\scrO$ in place of $A$,
there are $ x_U \iin G $ for $ U \iin \scrO $ such that
$ \psm ( F_U x_U , F_V x_V ) \!>\! 1 $ whenever
$ U, V \iin \scrO $, $ U \neq V $.
Define the function $ f : G \rightarrow \real $ by
\[
f(x) = \sup_{V \in \scrO} \; \max_{y \in F_V} \;
( \; h_V (y) - \psm( x , y x_V ) \; )^+
\;\;\; \mbox{\rm for} \;\;\; x \in G.
\]
Each function $ \sect{x}{( \; h_V (y) - \psm( x , y x_V ) \; )^+} $ belongs
to $ \Lipplus(\psm) $, and thus $ f \in \Lipplus(\psm)$.
It remains to be proved that $ f (x x_U ) = h_U (x) $ for every $ U \iin \scrO $
and $ x \iin F_U $.
Once that is established, it will follow that
$ \rtrans{x_U}{f} \iin U $ for every $ U \iin \scrO $.

Take any $ U \iin \scrO $ and $ x \iin F_U $.
From the definition of $f$ we get $ f (x x_U ) \geq h_U (x) $.
To prove the opposite inequality, consider any $ V \iin \scrO $ and any $ y \iin F_V $.

Case I: $ V = U $.
From
$
h_U (y) - h_U (x) \leq | h_U (y) - h_U (x) | \leq \psm(x,y) \leq \psm( x x_U, y x_U )
$
and $ h_U = h_V $, $ x_U = x_V $, we get
$ ( h_V (y) - \psm (x x_U , y x_V ) )^+ \leq h_U (x) $.

Case II: $ V \neq U $.
From $ \psm (x x_U , y x_V ) > 1 $ we get
$ ( h_V (y) - \psm (x x_U , y x_V ) )^+ = 0 \leq h_U (x) $.

Thus
$ ( h_V (y) - \psm (x x_U , y x_V ) )^+ \leq h_U (x) $ in both cases,
and now $ f (x x_U ) \leq h_U (x) $ follows from the definition of $f$.
\qed

\begin{lemma}
    \label{lemma:density}
Let $G$ be a topological group and $ \psm \in \rp(G)$.
If $ d( \psm ) = \bccard(\psm) \!\geq \aleph_0 $
then there exists $ f \iin \Lipplus(\psm) $
such that $ \Lipplus(\psm) = \clorbit{f} $.
\end{lemma}

\noindent
{\bf Proof.}
Let $H$ be a $\psm$-dense subset of $G$ such that $ \card{H} = d ( \psm ) = \bccard(\psm) $.
On $\Lipplus(\psm)$, the topology of pointwise convergence on $G$ coincides with
the topology of pointwise convergence on $H$.
Thus $\Lipplus(\psm)$ is homeomorphic to a subset of the product space $ \real^H $,
its topology has a base of cardinality at most $ \bccard(\psm) $,
and by Lemma~\ref{lemma:largeorbit} there is $ f \iin \Lipplus(\psm) $
whose orbit intersects every nonempty open set in $ \Lipplus(\psm) $.
\qed

In this paper,
Lemma~\ref{lemma:density} is the key for finding sufficient conditions for ambitability.
If for every $\psm\iin\rp(G)$ there exists $\psm'\iin\rp(G)$ such that $\psm'\geq\psm$ and
$d(\psm')=\bccard(\psm')\!\geq\aleph_0$
then the group $G$ is ambitable.

\begin{theorem}
    \label{theorem:mainth}
Let $\kappa$ be an infinite cardinal,
and $G$ a locally $\kappa$-bounded topological group.
If there exists $ \psm_0 \iin \rp(G) $ such that $ \bcard(\psm_0) \geq \kappa $
then $G$ is ambitable.
\end{theorem}

\noindent
{\bf Proof.}
Take any $ \psm \iin \rp(G) $.
Let $e$ be the identity element of $G$, and
$ B = \{ x \iin G \; | \; \psm(e,x) \leq 1 \} $.
Without loss of generality, assume that
$ \bcard (\psm) \geq \kappa $
and $B$ has a $\psm$-dense subset $H$ such that $ \card{H} \leq \kappa $.
(If $\psm$ does not have these properties then replace $\psm$ by
a larger pseudometric in $ \rp (G) $ that does.)

By Lemma~\ref{lemma:cardftions}, there is $ P \subseteq G $ such that
$ \card{P} = \bcard (\psm) $ and $ G = BP $.
The set $HP$ is $\psm$-dense in $G$,
therefore $ d( 2 \psm ) = d(\psm) \leq \, \kappa \cdot \, \bcard (\psm) = \bcard (\psm)$.
Thus
$ d( 2 \psm ) \leq \bccard (2 \psm)$
by Theorem~\ref{theorem:bctobcc},
and by Lemma~\ref{lemma:density}
there is $ f \iin \Lipplus(2 \psm) $ such that $ \Lipplus(2 \psm) = \clorbit{f} $.
\qed

\begin{corollary}
    \label{corollary:succard}
Let $\kappa$ be an infinite cardinal.
Every topological group that is locally $\kappa^+$-bounded and not $\kappa$-bounded
is ambitable.
\end{corollary}

\noindent
{\bf Proof.}
Let $G$ be locally $\kappa^+$-bounded and not $\kappa$-bounded.
By Lemma~\ref{lemma:bounded}
there is $\psm \iin \rp(G) $ for which $\bcard(\psm) \geq \kappa^+$,
and Theorem~\ref{theorem:mainth} applies with $\kappa^+$ in place of $\kappa$.
\qed

\begin{theorem}
    \label{theorem:lanambitable}
When $n$ is a positive integer,
every locally $\aleph_n$-bounded topological group is either precompact or ambitable.
\end{theorem}

\noindent
{\bf Proof.}
Let $G$ be locally $\aleph_n$-bounded for some $n$.
Let $m \geq 0$ be the smallest integer for which $G$ is locally $\aleph_m$-bounded.
If $m \geq 1$ then $G$ is ambitable by Corollary~\ref{corollary:succard}
with $\kappa = \aleph_{m-1}$.
If~$m = 0$ and $G$ is not precompact
then there exists $ \psm \iin \rp(G) $ such that $ \bcard(\psm) \geq \aleph_0 $
and $G$ is ambitable by Theorem~\ref{theorem:mainth}.
\qed

\begin{corollary}
    \label{corollary:lcambitable}
Every locally compact topological group is either compact or ambitable.
\end{corollary}

\begin{corollary}
    \label{corollary:obambitable}
Every $\aleph_0$-bounded topological group is either precompact or ambitable.
\end{corollary}

Lemma~\ref{lemma:density} yields also other classes of ambitable groups,
such as those in the next two theorems.

\begin{theorem}
    \label{theorem:additivegroup}
If $G$ is the additive group of an infinite normed space with the norm topology
then $G$ is ambitable.
\end{theorem}

\noindent
{\bf Proof.}
Let $\|\cdot\|$ be the norm, and $\psm(x,y)=\|x-y\|$, $x,y\in G$.
The topology of $G$ is defined by the metric $\psm$.
Since $ \bccard(\psm) =  \bccard(k \psm) \geq \aleph_0 $ for $ k=1,2, \ldots $,
it follows that $d(\psm)=\bccard(\psm)$ by Lemma~\ref{lemma:cardftions}.

Let $\psm'$ be any pseudometric in $\rp(G)$.
Without loss of generality, assume that $\psm'\geq\psm$
(if not then replace $\psm'$ by $\psm+\psm'$).
Since
$\bccard(\psm')\geq\bccard(\psm)=d(\psm)\geq d(\psm')$,
by Lemma~\ref{lemma:density} there exists
$ f \iin \Lipplus(\psm' ) $
such that $ \Lipplus(\psm' ) = \clorbit{f} $.
\qed

Let $\kappa$ be an infinite cardinal.
Define $\cf (\kappa)$, the {\em cofinality\/} of $\kappa$,
to be the smallest cardinality of a set $A$ of cardinals such that
$\kappa' < \kappa$ for each $ \kappa' \in A $ and $ \sup A = \kappa $.
Jech (\cite{Jech2002},~1.3) discusses cofinality in detail.

\begin{theorem}
    \label{theorem:forcofinality}
If $G$ is a topological group and
for every $\psm\iin\rp(G)$ there exists $\psm'\iin\rp(G)$
such that $\psm'\geq\psm$ and $\cf ( d(\psm') ) > \aleph_0 $,
then $G$ is ambitable.
\end{theorem}

\noindent
{\bf Proof.}
Let $\psm' \iin \rp(G) $ be such that $\cf ( d(\psm') ) > \aleph_0 $.
Then $ d(\psm') > \aleph_0 $,
and therefore
$ d(\psm') = \lim_{k \rightarrow \infty} \bccard(k \psm') $
by Corollary~\ref{corollary:bctobcc}.
Since $\cf ( d(\psm') ) > \aleph_0 $, it follows that
$ d(\psm') = \bccard(k \psm') $ for some $k$.
Thus $ d(k \psm' ) = \bccard(k \psm') $
and by Lemma~\ref{lemma:density} there exists
$ f \iin \Lipplus(k \psm') $
such that $ \Lipplus(k \psm') = \clorbit{f} $.
\qed

Note that $\cf ( \aleph_\omega ) = \aleph_0 $.
It is an open question whether
every $\aleph_\omega$-bounded topological group
is either precompact or ambitable.


\section{Topological centres}
    \label{section:centres}

Recall~\cite{Pachl2006} that for a topological group $G$,
\begin{itemize}
\item
$\Measru{G}$ is the norm dual of $\Ubru{G}$;
\item
the subspace $\UMeasru{G}$ of $\Measru{G}$ is defined as follows:
$\mu \in \UMeasru{G}$ iff $\mu$ is continuous on $\Lipplus(\psm)$
for each $\psm \in \rp(G)$;
\item
for $ x \in G $, $ \delta_x \in \Measru{G} $ is defined by
$ \delta_x (f) = f(x) $ for $ f \in \Ubru{G} $;
\item
the mapping $ \delta: x \mapsto \delta_x $ is a topological embedding of $G$
to $ \Measru{G} $ with the \wstar topology;
\item
when $ \mu , \nu \in \Measru{G} $, the {\em convolution\/}
of $\mu$ and $\nu$ is defined
by $ \mu \conv \nu (f)  = \mu ( \sect{x} \nu ( \sect{y} f(xy ) ) ) $
for $f \in \Ubru{G}$;
\item
$ \compactnru{G} $, the {\em uniform semigroup compactification of $G$\/},
is the \wstar closure of $ \delta(G) $ in $ \Measru{G} $,
with the \wstar topology and the convolution operation $\conv$ ;
\item
$ \compln{\ru G} = \compactnru{G} \cap \UMeasru{G}$ is a completion of $\ru G$.
\end{itemize}

A semigroup $S$ with a topology is a {\em right topological semigroup} if the mapping
$ \sect{x} xy $ from $S$ to $S$ is continuous for each $ y \in S $ (\cite{BJM1989}, 1.3).
For any right topological semigroup $S$, define its {\em topological centre\/}
\[
\Lambda (S) =  \{ \; x \in S \;
| \;
\mbox{\rm the mapping} \;\; \sect{y} xy  \;\;
\mbox{\rm is continuous on} \; S \; \} \; .
\]

The spaces studied here are denoted by various symbols in the literature.
Some of the more common notations are:
\begin{center}
\begin{tabular}[b]{|l|l|} \hline
{in this paper} \hspace{0.5cm} & {alternative notations} \\ \hline
$ \Ubru{G} $ & LUC$(G)$ or $ RUC\raisebox{1mm}{$\ast$} (G)$
  or $ {\cal U}_r (G) $ or $ {\cal LC} (G) $ \\ \hline
$ \Measru{G} $ & LUC$(G)\raisebox{1mm}{$\ast$} $
  or $ {\cal U}_r (G)\raisebox{1mm}{$\ast$} $
  or $ {\cal LC} (G)\raisebox{1mm}{$\ast$} $\\ \hline
$ \UMeasru{G} $  &  Leb$(G)$  \\ \hline
$ \compactnru{G} $  & $ G^{\mbox{\rm\scriptsize LUC}} $ or $ G^{\cal LC} $ \\ \hline
$ \Lambda ( \Measru{G} ) $  & $ Z(G) $ or $ Z_t (G) $ \\ \hline
\end{tabular}
\end{center}

Let $G$ be a topological group.
Then $ \Measru{G} $ with the $\conv\;$ operation and
the \wstar topology is a right topological semigroup.
This semigroup and its subsemigroup $ \compactnru{G}$
have an important role in functional analysis on $G$.
Significant research efforts have been devoted to characterizing
their topological centres.
In the rest of this section we will see how the known results follow from results about
ambitable groups.
The same approach yields a characterization of topological centres not only
in $\Measru{G}$ and $\compactnru{G}$ but also in any intermediate semigroup between
$\Measru{G}$ and $\compactnru{G}$.

Research in abstract harmonic analysis is often concerned with linear spaces
over $\mathbb{C}$, the field of complex numbers, rather than the field $\real$
used here.
However, it is an easy exercise to derive the $\mathbb{C}$-version of any result in
this paper from its $\real$-version.

For every topological group $G$ we have $ \UMeasru{G} \subseteq
\Lambda ( \Measru{G} ) $; see Proposition~4.2
in~\cite{Ferri-Neufang2007} or section~5 in~\cite{Pachl2006}. If
$G$ is precompact then $ \UMeasru{G} = \Measru{G} $ and therefore
$ \UMeasru{G} = \Lambda ( \Measru{G} ) $.

\begin{question}
    \label{question:centreM}
Is $ \UMeasru{G} = \Lambda ( \Measru{G} ) $ for every topological group $G$ ?
\end{question}

The positive answer was proved for locally compact groups by Lau~\cite{Lau1986}
and for $\aleph_0$-bounded groups by
Ferri and Neufang~\cite{Ferri-Neufang2007}.
By Corollary~\ref{corollary:twosemigroups} below,
the answer is positive for every ambitable group.

The situation is similar for $ \Lambda ( \compactnru{G} ) $.
For every topological group $G$ we have
$ \compln{ \ru G } \subseteq \Lambda ( \compactnru{G} ) $.
If $G$ is precompact then $ \compln{ \ru G } = \compactnru{G} $
is a compact group and therefore $ \compln{ \ru G } = \Lambda ( \compactnru{G} ) $.

\begin{question}
    \label{question:centreC}
Is $ \compln{ \ru G } = \Lambda ( \compactnru{G} ) $ for every topological group $G$ ?
\end{question}

The positive answer was proved for locally compact groups by Lau and Pym~\cite{Lau-Pym1995}
and for $\aleph_0$-bounded groups by
Ferri and Neufang~\cite{Ferri-Neufang2007}.
Again the answer is positive for every ambitable group, by
Corollary~\ref{corollary:twosemigroups}.

\begin{lemma}
    \label{lemma:halfconv}
Let $G$ be a topological group and $ f \in \Ubru{G}) $.
\begin{enumerate}
\item
The mapping $ \varphi : \nu \mapsto \sect{x} \nu ( \sect{y} f (xy) ) $ is continuous
from $ \compactnru{G} $  to the product space $ \real^G $.
\item
$ \varphi ( \compactnru{G} ) = \clorbit{f} $.
\end{enumerate}
\end{lemma}

\noindent
{\bf Proof.}
1. As noted above, $ \delta_x \in \Lambda ( \compactnru{G} ) $ for each $ x \in G $,
and thus the mapping
$ \nu \mapsto \delta_x \conv \nu $ is \wstar continuous from
$ \compactnru{G} $ to itself.
Since $ \delta_x \conv \nu (f)  = \nu ( \sect{y} f(xy) )$,
this means that the mapping
$ \nu \mapsto \nu ( \sect{y} f(xy) ) $ from
$ \compactnru{G} $ to $\real$ is continuous for each $ x \in G $,
and therefore
the mapping $ \nu \mapsto \sect{x} \nu ( \sect{y} f (xy) ) $ is continuous
from $ \compactnru{G} $ to $ \real^G $.

2. $ \varphi ( \delta_x ) = \rtrans{x}{f} $ for all $ x \in G $,
and therefore $ \varphi ( \delta ( G ) ) = \orbit{f} $.
The mapping $ \varphi $ is continuous by part 1, $ \compactnru{G} $ is compact, and
$ \delta(G) $ is dense in $ \compactnru{G} $.
It follows that $ \varphi ( \compactnru{G} ) = \clorbit{f} $.
\qed

\begin{lemma}
    \label{lemma:orbitcontinuity}
Let $G$ be any topological group, $ \mu \in \Measru{G} $ and $ f \in \Ubru{G} $.
If the mapping $ \nu \mapsto \mu \conv \nu $ from $ \compactnru{G} $ to $ \Measru{G} $
is \wstar continuous then $ \mu $ is continuous on $ \clorbit{f} $.
\end{lemma}

\noindent
{\bf Proof.}
As in Lemma~\ref{lemma:halfconv}, define
$ \varphi ( \nu ) =  \sect{x} \nu ( \sect{y} f (xy) ) $
for $ \nu \in \compactnru{G} $.

\begin{picture}(300,90)
\thicklines
\put(50,60){\makebox(30,20){$\compactnru{G}$}}
\put(140,60){\makebox(40,20){$\clorbit{f}$}}
\put(140,0){\makebox(40,20){$\real$}}

\put(80,70){\vector(1,0){60}}
\put(160,60){\vector(0,-1){40}}
\put(80,60){\vector(3,-2){66}}

\put(95,66){\makebox(30,20){$\varphi$}}
\put(160,30){\makebox(20,20){$\mu$}}
\put(50,20){\makebox(80,20){$\nu \mapsto \mu \conv \nu (f) $}}
\end{picture}

By the definition of convolution,
$ \mu \conv \nu (f)  = \mu ( \sect{x} \nu ( \sect{y} f(xy ) ) )
= \mu ( \varphi ( \nu ) ) $.
Thus $ \mu \circ \varphi $ is continuous from $ \compactnru{G} $ to $\real$.

By Lemma~\ref{lemma:halfconv},
$ \varphi $ is continuous
from $ \compactnru{G} $ to $ \clorbit{f} $,
and $ \varphi ( \compactnru{G} ) = \clorbit{f} $.
Since $ \compactnru{G} $ is compact, it follows that
$\mu$ is continuous on $\clorbit{f}$.
\qed

\begin{theorem}
    \label{theorem:ambitableunif}
If $G$ is an ambitable topological group,
$ S \subseteq \Measru{G} $,
and $S$ with the $\conv$ operation is a semigroup
such that $ \compactnru{G} \subseteq S $,
then $ \Lambda (S) = \UMeasru{G} \cap S $.
\end{theorem}

\noindent
{\bf Proof.}
As was noted above, $ \UMeasru{G} \subseteq \Lambda ( \Measru{G} $.
Therefore $ \UMeasru{G} \cap S \subseteq \Lambda (S) $ for every semigroup
$ S \subseteq \Measru{G} $.

To prove the opposite inclusion,
take any $ \mu \iin \Lambda (S) $ and any $ \psm \iin \rp(G) $.
Since $ \compactnru{G} \subseteq S $,
the mapping $ \nu \mapsto \mu \conv \nu $ from $ \compactnru{G} $ to $ \Measru{G} $
is \wstar continuous by the definition of $ \Lambda (S) $.
Since $ G $ is ambitable, $ \Lipplus(\psm) \subseteq \clorbit{f} $
for some $ f \in \Ubru{G} $.
By~\ref{lemma:orbitcontinuity}, $\mu$ is continuous on  $ \clorbit{f} $ and therefore
also on $ \Lipplus(\psm) $.
Thus $ \mu \in \UMeasru{G} $.
\qed

\begin{corollary}
    \label{corollary:twosemigroups}
If $G$ is an ambitable topological group then
$ \UMeasru{G} = \Lambda ( \Measru{G} ) $
and
$ \compln{ \ru G } = \Lambda ( \compactnru{G} ) $.
\end{corollary}

\noindent
{\bf Proof.}
Apply \ref{theorem:ambitableunif} with $ S = \Measru{G} $ and with $ S = \compactnru{G}$.
\qed

Combined with Theorem~\ref{theorem:lanambitable},
Corollary~\ref{corollary:twosemigroups} gives positive answers to
Questions~\ref{question:centreM} and~\ref{question:centreC} for
locally $\aleph_0$-bounded groups,
and in particular for locally compact and for $\aleph_0$-bounded groups.
This generalizes the previously published results cited above.

Note that Theorem~\ref{theorem:ambitableunif} applies not only
to the semigroups $ \compactnru{G} $ and $ \Measru{G} $ in
Corollary~\ref{corollary:twosemigroups}, but also to many other
semigroups between $ \compactnru{G} $ and $ \Measru{G} $
--- for example, the semigroup of all positive elements in $ \Measru{G} $,
or the semigroup of all finite linear combinations of elements of $ \compactnru{G} $
with integral coefficients.

\begin{corollary}
    \label{corollary:uamenable}
No uniquely amenable topological group is ambitable.
\end{corollary}

\noindent
{\bf Proof.}
If a group $G$ is ambitable then $ \UMeasru{G} = \Lambda ( \Measru{G} ) $
by Corollary~\ref{corollary:twosemigroups}.
If $G$ were also uniquely amenable, then $G$ would be precompact by
Theorem~5.2 in~\cite{Pachl2006},
which would contradict Theorem~\ref{theorem:precompact} above.
\qed

\noindent
{\bf Acknowledgement.\/}
I wish to thank Matthias Neufang and Vladimir Uspenskij for beneficial discussions
of the concepts presented in this paper,
and Juris Stepr\={a}ns for sharing a helpful counterexample.


\end{document}